# RNGA for non-square multivariable control systems: properties and application

Shaival Hemant Nagarsheth* and Shambhu Nath Sharma

*Electrical Engineering Department, Sardar Vallabhbhai National Institute of Technology, Surat, Gujarat 395007, India.*
*E-mail address: shn411@gmail.com, snsvolterra@gmail.com.*

**Abstract**

The Relative Gain Array (RGA) and Relative Normalized Gain Array (RNGA) have received considerable attention for square systems. In this paper RNGA $\lambda^{RN}$ with the column major, i.e. $\lambda^{RN} \in R^{r \times s}, r < s,$ for non-square multivariable systems is introduced. RNGA of the paper has a row-column inequality, i.e. $r$ number of rows is less than the $s$ number of columns. Unlike the conventional RGA, the RNGA loop pairing criteria of the paper considers both steady-state as well as transient information for the assessment of control-loop interactions. The RNGA for square systems is extended for non-square multivariable systems by thoroughly deriving its supporting properties. The RNGA method is applied to a non-square multivariable radiator laboratory test setup for loop pairing. Closed-loop results arising from the RNGA-based loop pairing are depicted in the paper. The lacuna of the conventional RGA loop pairing has been overcome by the application of the developed RNGA of this paper. The results unfold the effectiveness of RNGA over RGA for non-square multivariable systems to have minimum interactions and better control.

## Introduction

For multi-input and multi-output systems, where the number of outputs and number of inputs is unequal, non-square transfer function matrices arise to be an adequate representation. The control-loop pairing is a potential problem for non-square systems to achieve minimum interactions among the control loops. The solution to the above problem lies with the decentralized control of the non-square multivariable system. The importance of a decentralizing control scheme lies with its simplicity to implement for field engineers. Many of the advanced control methods though better in performance are not feasible to implement due to its complex nature, long implementation time, lack of understanding [1]. Hence, advancement in the decentralized control theory warrants investigation. The primary step in the decentralized control is to decide the control-loop pairing based on the loop-interactions measure. The suitable pairing based on the minimal loop-interactions confirms that the closed-loop performance is not deteriorated [2]. Bristol [3] introduced a relative gain array notion to decide the loop pairing for square systems. The method chiefly involves the construction of the gain array symmetric matrix using the Schur product. Further, the RGA properties were sketched by Grosdidier et al. [4]. Significantly, the RGA properties for non-square systems with more outputs and less inputs were derived by Chang and Yu [5].

Cameron [6] presented interaction indicators for square multivariable systems with the knowledge of eigenstructure. The control-loop pairing decision through RGA leads to decentralized control of the system [7].

Generally, the transient information of the system is not accounted for attaining the gain matrix. McAvoy et al. [8] presented a new approach of Dynamic Relative Gain Array (DRGA) to give more accurate interactions amongst loops. The Relative Disturbance Gain Array (RDGA) was introduced to account for disturbance effects [9]. The universality of the RGA is not restricted to control-loop pairing only, RGA coupled with condition number has proven useful for robustness analysis of multivariable systems [10,11]. The Relative Effective Gain Array (REGA) [12] does not depend on controller parameters but on the critical frequency of the individual element. This is regarded as the Relative Gain Array from the frequency perspective. He et al. [13] presented a RNGA based new control-loop pairing criterion for square multivariable systems. The RNGA accounts for steady-state gain, dead time and time constant parameters. RNGA has been beneficial in comparison to the conventional RGA, as it takes transient information into account. As compared to DRGA, RNGA provides loop interactions independent of the controller type. Moreover, the computation of RNGA has proven simple and useful to work with, for the field engineers to carry out the loop pairing decision for decentralized control of practical industrial problems [13]. RNGA can also be a better predictor of system interactions for non-square practical systems. The universality of the non-square system with the less output-more input structure can be found in practical problems: mixing-tank process (2 outputs and 3 inputs) [14], Air path scheme of a turbocharged diesel engine (2 outputs and 3 inputs) [15,16], hot oil fractionators (2 outputs and 4 inputs) [17], Shell control problem (5 outputs 7 inputs) [18] etc. Hence, defining the RNGA for non-square multivariable systems with the rigorous derivation of its supporting properties can become a worthy contribution.

This paper presents the RNGA for non-square multivariable systems with application to a non-square radiator control system having less output-more input in its setup. A systematic derivation of RNGA properties is carried out in the paper such that it can be applied to a non-square multivariable radiator system. The minimal pairing achieved here pinpoints to a quick and superior control setting i.e. achieving the desired, through fewer control efforts. The RNGA property for the $r < s$ case, where $r$ denotes the output dimension and $s$ denotes the input dimension, suggests the input variables elimination from control configuration selections. The method adopted in this paper unifies Schur product algebra [19], generalized inverse [19] and the calculus of matrices to weave the non-square RNGA properties, since the RNGA is a non-square matrix with the column major for the non-square system with the less output-more input structure. To test the effectiveness, the technique is applied to a non-square multivariable radiator laboratory setup. The system transfer matrix is derived from the experimental step test readings for the radiator laboratory setup. The control-loop pairing of the non-square system is carried using the theory developed in this paper. IMC (Internal Model Control) tuned PID controllers are designed for the decentralized control of the non-square system. The closed-loop performance resulting from the suggested configuration by non-square RNGA is compared to that of the control configuration of the conventional RGA.

*Notations*: Consider $G(s)$ is a $r \times s$ non-square process transfer matrix with $r < s$, where $Y(s) = G(s)U(s)$. It is important to note that throughout the paper, we weave the RNGA properties for the case $r < s$. The term column major of the paper indicates that, the non-square matrix has more columns with less associated rows, i.e. more inputs and less outputs.

**Definition of the non-square RNGA**

Suppose $Y(s)$ is an $r \times 1$ output vector and $U(s)$ is an $s \times 1$ input vector. Now, $G(s) = (G_{ij}) = \left( \dfrac{k_{ij} e^{-(t_d)_{ij} s}}{1 + (\tau_c)_{ij} s} \right)$, where $k_{ij}, (t_d)_{ij}, (\tau_c)_{ij}$ are the process gain, dead time and time constant respectively of the $i^{th}$ output with respect to the $j^{th}$ input and the variables, $i$ and $j$, run over 1 to $r$ and 1 to $s$. The steady state gain and the average residence time [20] parameter of the transfer function entries of the non-square transfer matrix are

$$k = (k_{ij}) = (G_{ij}(0)), \quad b_{ij} = (\tau_c)_{ij} + (t_d)_{ij}. \tag{1a}$$

*Remark 1:* Consider Second Order Plus Dead Time (SOPDT) system $G(s) = (G_{ij}) = \left( \dfrac{k_{ij} e^{-(t_d)_{ij} s}}{(1+(\tau_1)_{ij} s)(1+(\tau_2)_{ij} s)} \right)$, where $(\tau_1 \neq \tau_2)$. The average residence time will be $b_{ij} = (\tau_1 + \tau_2)_{ij} + (t_d)_{ij} = (\tau_c)_{ij} + (t_d)_{ij}$.

Now, introducing an alternative and a computational convenient notation, $\tilde{B} = (\tilde{B}_{ij}) = (b_{ij}^{-1})$, the Normalized Gain Array (NGA) can be defined as $A = K \circ \tilde{B}$, where $\circ$ is the Schur product between two non-square matrices, where the size of their matrices coincides. Thus,

$$A_{ij} = k_{ij} \tilde{B}_{ij} = k_{ij} b_{ij}^{-1} = \dfrac{k_{ij}}{(\tau_c)_{ij} + (t_d)_{ij}}. \tag{1b}$$

Using generalized inverse result ([21], p. 523), the non-square RNGA becomes

$$\lambda^{RN} = A \circ (A^+)^T = A \circ (A^T(AA^T)^{-1})^T = A \circ (AA^T)^{-1} A. \tag{2}$$

On the other hand, non-square RNGA can be represented in terms of system parameters. In the component-wise setting, equation (2) can be recast as

$$\lambda^{RN} = (\lambda_{ij}^{RN}) = (A_{ij}((AA^T)^{-1} A)_{ij}). \tag{3}$$

After combining equation (3) with equation (1), we get

$$\lambda_{ij}^{RN} = A_{ij}\sum_{\gamma}((\sum_{\phi}A_{i\phi}A_{j\phi})^{-1})_{i\gamma}A_{\gamma j}$$

$$= \frac{k_{ij}}{(\tau_c)_{ij}+(t_d)_{ij}}\sum_{\gamma}\left(\left(\sum_{\phi}\frac{k_{i\phi}}{(\tau_c)_{i\phi}+(t_d)_{i\phi}}\frac{k_{j\phi}}{(\tau_c)_{j\phi}+(t_d)_{j\phi}}\right)^{-1}\right)_{i\gamma}\frac{k_{\gamma j}}{(\tau_c)_{\gamma j}+(t_d)_{\gamma j}}.$$

$$= \sum_{\gamma}\left(\left(\sum_{\phi}\frac{k_{i\phi}}{(\tau_c)_{i\phi}+(t_d)_{i\phi}}\frac{k_{j\phi}}{(\tau_c)_{j\phi}+(t_d)_{j\phi}}\right)^{-1}\right)_{i\gamma}\frac{k_{\gamma j}}{(\tau_c)_{\gamma j}+(t_d)_{\gamma j}}\frac{k_{ij}}{(\tau_c)_{ij}+(t_d)_{ij}}. \quad (4)$$

*Remark 2:* The non-square RNGA expression of equation (4) is valid for $r \times s$ transfer matrix $G(s)$, where $r < s$. To construct the non-square RNGA expression for the case $r > s$, $A^+ = (A^T A)^{-1}A^T$ ([22], p. 39) and $\lambda^{RN} = A \circ (A^+)^T = A \circ ((A^T A)^{-1}A^T)^T = A \circ A(A^T A)^{-1}$ are used. An alternative non-square RNGA expression for the case $r > s$, can be stated as

$$\lambda^{RN} = (\lambda_{ij}^{RN}) = (A_{ij}(A(A^T A)^{-1})_{ij}). \quad (5a)$$

Alternatively,

$$\lambda_{ij}^{RN} = A_{ij}\sum_{\gamma}((\sum_{\phi}A_{\phi i}A_{\phi j})^{-1})_{\gamma j}A_{i\gamma} = \sum_{\gamma}((\sum_{\phi}A_{\phi i}A_{\phi j})^{-1})_{\gamma j}A_{i\gamma}A_{ij}.$$

On combining equation (1) with the term $\lambda_{ij}^{RN}$, we have

$$\lambda_{ij}^{RN} = \sum_{\gamma}\left(\left(\sum_{\phi}\frac{k_{\phi i}}{(\tau_c)_{\phi i}+(t_d)_{\phi i}}\frac{k_{\phi j}}{(\tau_c)_{\phi j}+(t_d)_{\phi j}}\right)^{-1}\right)_{\gamma j}\frac{k_{i\gamma}}{(\tau_c)_{i\gamma}+(t_d)_{i\gamma}}\frac{k_{ij}}{(\tau_c)_{ij}+(t_d)_{ij}}. \quad (5b)$$

The above equation (5b) would be useful for the controller pairing of non-square systems with less input-more output configurations.

**Properties of the non-square RNGA**

Consider the RNGA $\lambda^{RN}$ contained in the real-valued matrix space $R^{r \times s}$, where $r < s$. The row sum $R(i)$, the column sum $C(j)$, scaling and permutation properties of the non-square RNGA are derived. Here, we list the following useful RNGA relations:

$$\lambda_{ij}^{RN} = A_{ij}((A^+)^T)_{ij} = A_{ij}A_{ji}^+, \quad (6)$$

$$\sum_{j}\lambda_{ij}^{RN} = \sum_{j}A_{ij}A_{ji}^+ = \sum_{\phi}A_{i\phi}A_{\phi i}^+ = \sum_{\substack{\phi \\ j=i}}A_{i\phi}A_{\phi j}^+ = (AA^+)_{ii} = R(i), \quad (7a)$$

$$\sum_i \lambda_{ij}^{RN} = \sum_i A_{ij} A_{ji}^+ = \sum_{\substack{\phi \\ i=j}} A_{\phi j} A_{j\phi}^+ = \sum_{\substack{\phi \\ i=j}} A_{\phi j} A_{i\phi}^+ = \sum_\phi A_{i\phi}^+ A_{\phi j} = (A^+ A)_{jj} = C(j). \qquad (7b)$$

**Property 1:** Consider the RNGA of an $r \times s$ matrix $A$ is $\lambda^{RN}$, where $r < s$. The sum of elements in each row of the non-square RNGA $\lambda^{RN}$ is always equal to '1', $R(i) = 1$, $1 \le i \le r$.

*Proof:* Making the use of the generalized inverse result of Graybill et al. ([21], p. 523), Penrose [23] and equation (7a), we get

$$R(i) = (AA^+)_{ii} = (A A^T (AA^T)^{-1})_{ii} = (I_r)_{ii} = 1.$$

**Property 2:** The sum of elements in each column of the non-square RNGA $\lambda^{RN}$ is always between zero and unity, i.e. $0 \le C(j) \le 1, 1 \le j \le s$. The RNGA $\lambda^{RN}$ is the $r \times s$, $r < s$.

*Proof:* Consider the non-square RNGA

$$\lambda_{ij}^{RN} = \frac{1}{2} \frac{A_{ij}}{\alpha} \frac{d\alpha}{dA_{ij}}, \qquad (8)$$

where $\alpha = \det(AA^T)$, the sizes of the matrices $AA^T$ and $A^T A$ are $r \times r$ and $s \times s$, where $r < s$. The Binet-Cauchy relation [24] for the case $r < s$, can be recast. Furthermore

$$\det(AA^T) = \sum_l \det(A^T A)_l = \sum_{1 \le \kappa \le \binom{s}{r}} (\det A(\kappa))^2,$$

$$A_{ij} = k_{ij} \tilde{B}_{ij} = k_{ij} b_{ij}^{-1} = \frac{k_{ij}}{(\tau_c)_{ij} + (t_d)_{ij}}.$$

Note that $1 \le \kappa \le \binom{s}{r}$, $l$ is a subset contained in the product space, where $1 \le l \le \binom{s}{r}$, $\kappa$ has $r$ tuples [25]. By using the matrix calculus ([26], p.169),

$$\frac{d\alpha}{dA_{ij}} = \frac{d \det(AA^T)}{dA_{ij}} = \frac{d}{dA_{ij}} \left( \sum_{1 \le \kappa \le \binom{s}{r}} (\det A(\kappa))^2 \right) = 2 \sum_{1 \le \kappa' \le \binom{s}{r}} \det A(\kappa) \frac{d \det A(\kappa)}{dA_{ij}}. \qquad (9)$$

Note that

$$\frac{d \det A(\kappa)}{dA_{ij}} = \sum_m \frac{d \nabla_{mn}(\kappa') A_{mn}(\kappa)}{dA_{ij}} = \sum_m \left( \frac{d \nabla_{mn}(\kappa')}{dA_{ij}} A_{mn}(\kappa) + \nabla_{mn}(\kappa') \frac{dA_{mn}(\kappa)}{dA_{ij}} \right) = \nabla_{in}(\kappa'). \qquad (10)$$

The term $\nabla_{in}(\kappa')$ is the cofactor of the element $A_{in}(\kappa)$, where $1 \leq n \leq r, 1 \leq i \leq r$, $1 \leq \kappa' \leq \binom{s-1}{r-1}$ and $1 \leq \kappa \leq \binom{s}{r}$ that must include the $j^{th}$ column of $(A_{ij})$. The size of $(\nabla_{in}(\kappa'))$ is $(r-1) \times (r-1)$. After combining equations (9) and (10), we get

$$\frac{d\alpha}{dA_{ij}} = 2 \sum_{1 \leq \kappa' \leq \binom{s-1}{r-1}} \det A(\kappa') \nabla_{in}(\kappa'). \tag{11a}$$

For the convenient notations, replace $\nabla_{in}(\kappa')$ with $\det(A^{ij}(\kappa'))$ in equation (11a). The matrix $A^{ij}(\kappa')$ has the size $(r-1) \times (r-1)$, which is a consequence of $\frac{dA_{mn}(\kappa)}{dA_{ij}}$ and $A(\kappa')$ is an $r \times r$ matrix. Thus,

$$\frac{d\alpha}{dA_{ij}} = 2 \sum_{1 \leq \kappa' \leq \binom{s-1}{r-1}} \det A(\kappa') \nabla_{in}(\kappa') = 2 \sum_{1 \leq \kappa' \leq \binom{s-1}{r-1}} \det A(\kappa') \det(A^{ij}(\kappa')). \tag{11b}$$

Equation (7b) in conjunction with equations (8) and (11b) comes down to

$$\sum_i \lambda_{ij}^{RN} = \sum_i \frac{1}{2} \frac{A_{ij}}{\alpha} \frac{d\alpha}{dA_{ij}} = \frac{\sum_i A_{ij} \left( \sum_{1 \leq \kappa' \leq \binom{s-1}{r-1}} \det A(\kappa') \det(A^{ij}(\kappa')) \right)}{\sum_{1 \leq \kappa \leq \binom{s}{r}} (\det A(\kappa))^2}. \tag{12}$$

On further simplifications, equation (12) reduces to

$$\sum_i \lambda_{ij}^{RN} = \frac{\sum_{1 \leq \kappa' \leq \binom{s-1}{r-1}} \det A(\kappa') \sum_i A_{ij} \det(A^{ij}(\kappa'))}{\sum_{1 \leq \kappa \leq \binom{s}{r}} (\det A(\kappa))^2} = \frac{\sum_{1 \leq \kappa' \leq \binom{s-1}{r-1}} (\det A(\kappa'))^2}{\sum_{1 \leq \kappa \leq \binom{s}{r}} (\det A(\kappa))^2}. \tag{13}$$

The above equation (13) is a consequence of the Laplace expansion formula for the determinant. $\kappa$ runs over 1 to $\binom{s}{r}$ and $\kappa'$ runs over 1 to $\binom{s-1}{r-1}$. Each term within the summation sign is non-negative. Thus, the condition $0 \leq C(j) \leq 1$ holds.

**Property 3:** The non-square RNGA $\lambda^{RN}$ of an $r \times s$ matrix $A$ is output scaling-invariant, where $r < s$.

*Proof:* For the output scaling, an $r \times r$ diagonal matrix $Q_r$ is pre-multiplied to the normalized gain matrix $A$, where $Q_r = (q_i \delta_{ij}), 1 \leq i \leq r, 1 \leq j \leq r$. Since the property $(Q_r A)^+ = A^+ Q_r^{-1}$ holds [22], the output scaled RNGA can be written as

$$\overline{\lambda}^{RN} = (Q_r A) \circ (A^+ Q_r^{-1})^T = (Q_r A) \circ (Q_r^{-1}(A^+)^T).$$

Thus,

$$\overline{\lambda}_{ij}^{RN} = (\sum_\phi (Q_r)_{ii} \delta_{i\phi} A_{\phi j})(\sum_\phi (Q_r)_{ii}^{-1} \delta_{i\phi} A_{j\phi}^+) = (Q_r)_{ii} A_{ij} (Q_r)_{ii}^{-1} A_{ji}^+ = A_{ij} A_{ji}^+ = \lambda_{ij}^{RN}.$$

***Property 4:*** The non-square RNGA $\lambda^{RN}$ of an $r \times s$ matrix $A$ input scaling-variant, where $r < s$.

*Proof:* For the input scaling, an $s \times s$ diagonal matrix $Q_s$ is post-multiplied to the normalized gain matrix $A$. Suppose $Q_S = (q_i \delta_{ij}), 1 \leq i \leq s, 1 \leq j \leq s$. Here, we wish to calculate $(AQ_s)^+$. Making the use of the right generalized inverse $A^+ = A^T (AA^T)^{-1}$, we arrive at

$$(AQ_s)^+ = Q_s^T A^T (AQ_s^2 A^T)^{-1} = ((Q_s^T A^T (AQ_s^2 A^T)^{-1})^T)^T = (((AQ_s^2 A^T)^{-1})^T AQ_s)^T,$$

$$((AQ_s)^+)^T = ((AQ_s^2 A^T)^{-1})^T AQ_s = ((A^T)^+ (Q_s^2)^{-1} A^+)^T AQ_s$$

$$= (A^+)^T Q_s^{-1} Q_s^{-1} ((A^T)^+)^T AQ_s.$$

The $r \times r$ square matrix $(A^+)^T Q_s^{-1} Q_s^{-1} ((A^T)^+)^T$ can be rewritten as

$$(A^+)^T Q_s^{-1} Q_s^{-1} ((A^T)^+)^T = ((A^+)^T_{ij})(q_i^{-2} \delta_{ij})((A^T)^+_{ji}),$$

where $((A^+)^T_{ij})$ is an $r \times s$ matrix, $(q_i^{-2} \delta_{ij})$ is an $s \times s$ diagonal matrix, $((A^T)^+_{ji})$ is an $s \times r$ matrix. Thus, $(i, j)$ components of the $r \times r$ matrix $(A^+)^T Q_s^{-1} Q_s^{-1} ((A^T)^+)^T$ and the $r \times s$ matrix $AQ_s$ are

$$((A^+)^T Q_s^{-1} Q_s^{-1} ((A^T)^+)^T)_{ij} = \sum_{1 \leq \phi \leq s} \frac{A_{\phi i}^+ A_{\phi j}^+}{q_\phi^2}, \quad (AQ_s)_{ij} = A_{ij} q_j$$

respectively. Note that the notation $A_{ij}^+$ is the $(i, j)^{th}$ component of the right generalized inverse of the $r \times s$ matrix $A$, were $r < s$. Here, we wish to prove 'the non-square RNGA is input scaling-variant'. It suffices to prove that an element of the *input scaled* non-square RNGA $\hat{\lambda}^{RN}$ is different from that of the non-square RNGA $\lambda^{RN}$ at least. Thus, for the simplified notation, we show $\hat{\lambda}_{11}^{RN} \neq \lambda_{11}^{RN}$. The RNGAs are

$$\hat{\lambda}^{RN} = (AQ_s) \circ (A^+)^T Q_s^{-1} Q_s^{-1} ((A^T)^+)^T AQ_s, \quad \lambda^{RN} = A \circ (A^+)^T.$$

$$\hat{\lambda}^{RN}_{ij} = (AQ_s)_{ij} ((A^+)^T Q_s^{-1} Q_s^{-1} ((A^T)^+)^T AQ_s)_{ij}, \quad \lambda^{RN}_{ij} = A_{ij} A^+_{ji}.$$

Alternatively,

$$\hat{\lambda}^{RN}_{11} = (AQ_s)_{11} ((A^+)^T Q_s^{-1} Q_s^{-1} ((A^T)^+)^T AQ_s)_{11} = A_{11} \sum_{1 \le j \le r} \sum_{1 \le \phi \le s} \frac{A^+_{\phi 1} A^+_{\phi 1}}{q^2_\phi} A_{j1} q^2_1, \quad \lambda^{RN}_{11} = A_{11} A^+_{11}.$$

Thus, $\hat{\lambda}^{RN}_{11} \ne \lambda^{RN}_{11}$.

**Property 5:** Suppose $P_r$ and $P_s$ are two orthogonal matrices of the sizes $r \times r$ and $s \times s$ respectively. Consider the RNGA of an $r \times s$ matrix $A$ is $\lambda^{RN}$, and construct the matrix $P_r A P_s$, the RNGA $\tilde{\lambda}^{RN}$ of the matrix $P_r A P_s$ satisfies the following condition:

$$\tilde{\lambda}^{RN} = P_r \lambda^{RN} P_s,$$

where $\lambda^{RN} = A \circ (A^+)^T$.

*Proof:* The orthogonal matrices have properties, $P_s^{-1} = P_s^T$, $P_r^{-1} = P_r^T$. Now, we apply the definition of the non-square RNGA to the matrix $P_r A P_s$, we have

$$\tilde{\lambda}^{RN} = (P_r A P_s) \circ ((P_r A P_s)^+)^T = (P_r A P_s) \circ (P_s^+ A^+ P_r^+)^T = (P_r A P_s) \circ (P_s^{-1} A^+ P_r^{-1})^T.$$

After introducing the orthogonal matrix property and extending the property 2 of Grosdidier et al. [4], we are led to

$$\tilde{\lambda}^{RN} = (P_r A P_s) \circ (P_s^{-1} A^+ P_r^{-1})^T = (P_r A P_s) \circ (P_r (A^+)^T P_s) = P_r (A \circ (A^+)^T) P_s = P_r \lambda^{RN} P_s.$$

The above property reveals that the row and column permutations in the non-square normalized gain matrix $A$ introduce the same permutations in the non-square RNGA $\lambda^{RN}$. The property 5 is valid for both cases, $r < s$ and $r > s$.

**Application of RNGA to a non-square radiator control problem**

To demonstrate the usefulness of the non-square RNGA of the paper, a laboratory radiator setup is considered. First, the non-square transfer matrix of the radiator setup is achieved, then the control-loop pairing of the setup is adjudged on the basis of the non-square RNGA of the paper. Controller parameters are chosen using the IMC tuned PID controller for the radiator. The radiator schematic diagram is shown in Figure 1.

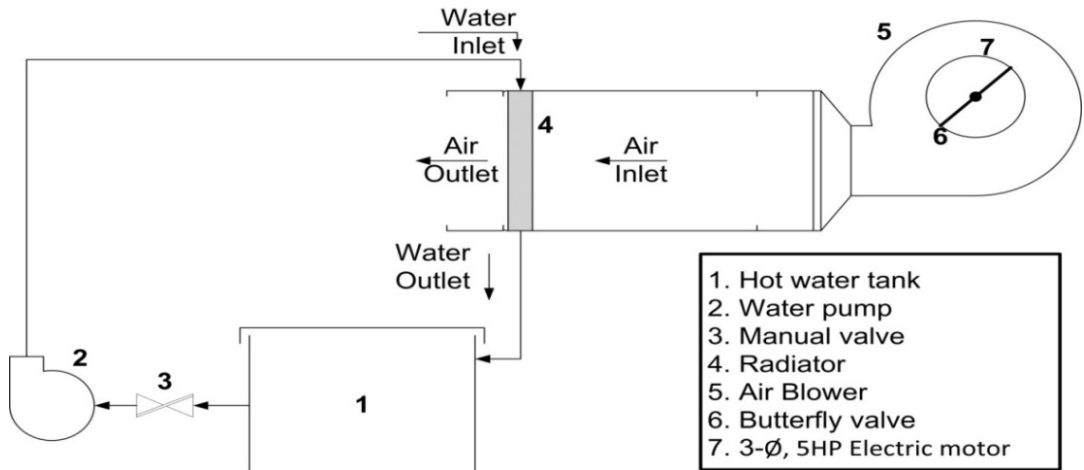
**Figure 1.** Schematic diagram of a laboratory radiator experimental setup

In the setup, water and air are fed from the hot water tank and the air blower housing is equipped with $3-\phi, 5\,\text{HP}$ electric motor. The heat transfer takes place between air and the hot water, which leads to cooling of the water passing through the radiator. At the radiator setup, four input variables are considered: air inlet temperature, air inlet flow rate, water inlet temperature, and water inlet flow rate. The output variables are air outlet and water outlet temperatures respectively. Figure 2 depicts the actual laboratory experimental setup. The system parameters of the radiator and operating points are given in Table 1.

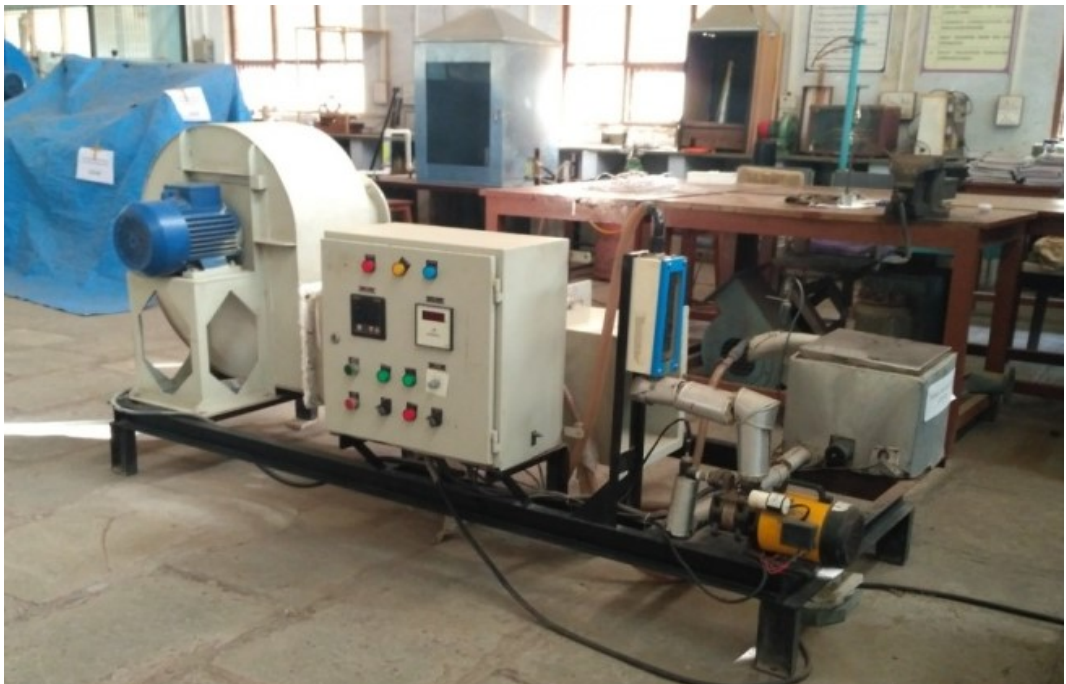
**Figure 2.** Laboratory radiator experimental setup

**Table 1.** Radiator setup parameters and operating point values

| Parameters | | Dimensions | Variables | Values |
|---|---|---|---|---|
| Radiator core | Height (mm) | 300 | Air inlet flow rate $(F_{AIN})\ m/s$ | 8.08 |
| | Length (mm) | 24 | Water inlet flow rate $(F_{WIN})\ LPM$ | 8 |
| | Width (mm) | 340 | Air inlet temperature $(T_{AIN})\ ^0C$ | 38.3 |
| Radiator fins (Aluminum) | Length (mm) | 24 | Water inlet temperature $(T_{WIN})\ ^0C$ | 72.5 |
| | Width (mm) | 10.89 | Air outlet temperature $(T_{AOUT})\ ^0C$ | 44.5 |
| | Thickness (mm) | 0.1 | Water outlet temperature $(T_{WOUT})\ ^0C$ | 65 |
| | Depth (mm) | 1.8 | | |
| | No. of fins | 185 | | |
| Radiator tube (Aluminum) | Outer diameter (mm) | 8 | | |
| | Inner diameter (mm) | 6 | | |
| | Thickness (mm) | 1 | | |

The step test method for the radiator non-square transfer matrix identification is adopted [27,28]. Each input is varied one by one considering other inputs at steady-state values. First, the air inlet flow rate $F_{AIN}$ is varied from its steady-state value, i.e. from $8.08\ m/s$ to $10\ m/s$, by keeping all other inputs at their steady-state. The effect of change in the air inlet flow rate to the two output temperatures $T_{AOUT}$ and $T_{WOUT}$ respectively is measured and noted to get the process reaction curve. Similarly, the water inlet flow rate is varied to get the relationship with the two output temperatures. Subsequently, experimental data from the open loop step test are taken for the rest of the two inputs. Using the procedure of the step test method, the open loop transfer matrix is arrived at, i.e.

$$(Y_1\ Y_2)^T = G(s)(U_1\ U_2\ U_3\ U_4)^T, \qquad (14)$$

where the size of the non-square transfer matrix $G(s)$ is $2 \times 4$.

$$G(s) = \begin{pmatrix} \dfrac{-0.9826e^{-13.74s}}{1+42.435s} & \dfrac{0.25702e^{-10.68s}}{1+32.922s} & \dfrac{1.09306e^{-18.67s}}{1+73.241s} & \dfrac{0.2154e^{-9.12s}}{1+78.7255s} \\ \dfrac{-0.1556e^{-7.971s}}{1+25.162s} & \dfrac{0.8045e^{-16.56s}}{1+30.264s} & \dfrac{0.3023e^{-19.86s}}{1+120.274s} & \dfrac{1.052e^{-18.27s}}{1+59.261s} \end{pmatrix},$$

Note that

$$(Y_1\ Y_2)^T = (T_{AOUT}\ T_{WOUT})^T,\ (U_1\ U_2\ U_3\ U_4)^T = (F_{AIN}\ F_{WIN}\ T_{AIN}\ T_{WIN})^T.$$

The closeness of the numerically simulated trajectories with the experimentally generated trajectories are depicted in Figure 3 and Figure 4. Figure 3 displays the process curve of output $Y_1(T_{AOUT})$ with respect to the change in all inputs. Figure 4 displays the same for output $Y_2(T_{WOUT})$. The model data is approximate and the experimental data reveal the exact physical situation. The difference is indicative of the radiator non-linearity, however, their non-linearity contribution is not appreciable. Thus, the radiator transfer matrix captures the qualitative characteristics of the non-square multivariable process considered here.

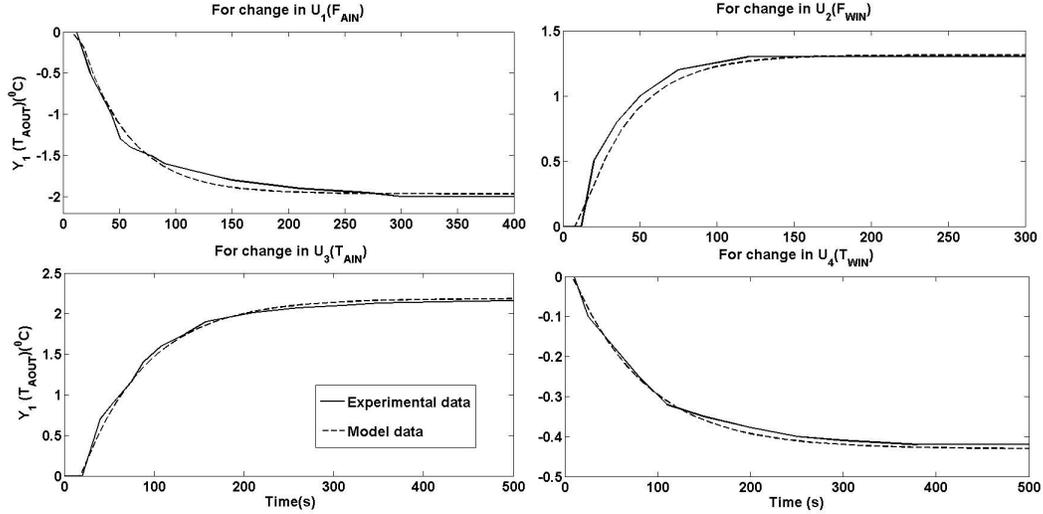

**Figure 3.** A Comparison of open loop step test for output 1, experimental and model

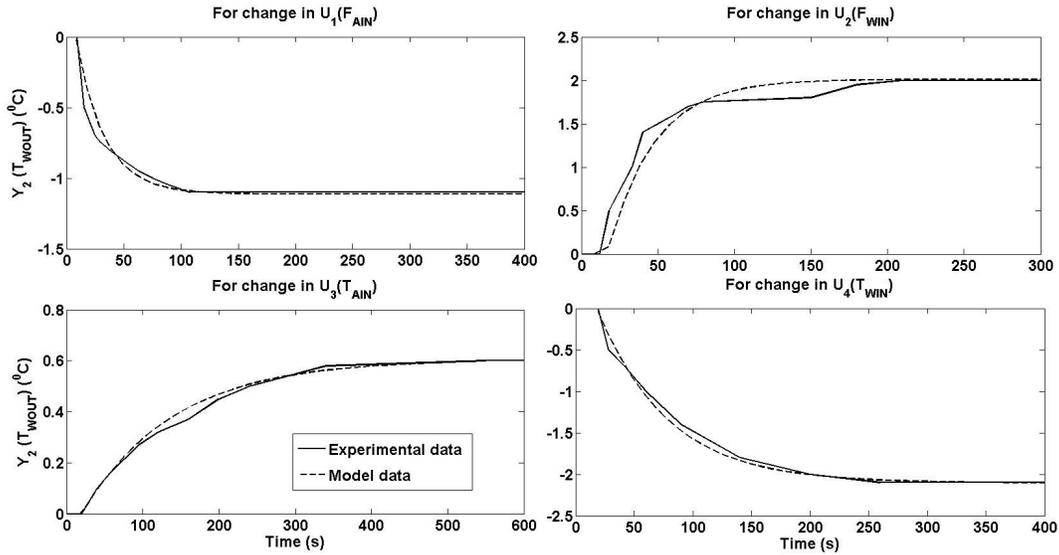

**Figure 4.** A Comparison of open loop step test for output 2, experimental and model

Now, the next step is to choose the best suitable control-loop pairing utilizing the non-square RNGA theory of the paper. That involves the following steps: (i) calculate the gain matrix, dead time matrix, time constant matrix and normalized gain matrix (ii) calculate non-square RNGA (iii) Consider the input-output pairing whose associated RNGA matrix element is $\geq 0.5$ and closer to unity [2].

Considering equations (1)-(4) in combination with the radiator non-square transfer matrix of equation (14), the resulting non-square radiator RNGA matrix is

$$\lambda^{RN} = \begin{pmatrix} 0.7166 & -0.0370 & 0.3470 & -0.0267 \\ -0.0486 & 0.6350 & -0.0210 & 0.4345 \end{pmatrix}. \tag{15}$$

The row sum property of the non-square radiator RNGA of equation (14) holds, i.e.

$$R(i) = (1.00 \quad 1.00)^T, \ 1 \leq i \leq 2. \tag{16a}$$

Furthermore, the column sum property of the non-square radiator RNGA of equation (15) also holds, i.e. $0 \leq C(j) \leq 1$,

$$(C(j))^T = (0.668 \ 0.598 \ 0.326 \ 0.4078)^T, 1 \leq j \leq 4, 0 \leq C(j) \leq 1. \tag{16b}$$

The input variant, output invariant properties as well as row-column permutations properties can be tested for the non-square radiator RNGA. Since the approach is direct, discussions are omitted. For the sake of comparisons, the non-square radiator RGA is discussed here. The non-square radiator RGA $\lambda^R$ is

$$\lambda^R = \begin{pmatrix} 0.4884 & -0.0194 & 0.5664 & -0.0354 \\ -0.0250 & 0.3759 & -0.0279 & 0.6770 \end{pmatrix}. \tag{17}$$

The non-square RGA is a consequence of equation (1a), equation (14) and the RGA definition $\lambda^R = k \circ (k^+)^T$. The row and column sum of the non-square radiator RGA of equation (16) are

$$R(i) = (1.00 \quad 1.00)^T, \tag{18a}$$

$$(C(j))^T = (0.4634 \ 0.3565 \ 0.5385 \ 0.6416)^T. \tag{18b}$$

The column sum values of the non-square RGA matrix and the value of $\lambda^R_{ij} \geq 0.5$ [2,29] recommend $(Y_1 - U_3 / Y_2 - U_4)$ pairing for the decentralized radiator control configuration. The first two inputs are eliminated from the control loop configuration, which is attributed to smaller entries of the column sum vector, see equation (18b). On the other hand, the non-square radiator RNGA of equation (15) recommends a different decentralized control scheme, i.e. the $(Y_1 - U_1 / Y_2 - U_2)$ pairing. Thus, the controller pairing $(Y_1 - U_1 / Y_2 - U_2)$ is an RNGA-based pairing. The controller pairing $(Y_1 - U_3 / Y_2 - U_4)$ is an RGA-based pairing.

To evaluate the closed-loop control performance of the non-square radiator system, decentralized IMC (Internal Model Control) tuned PID controllers are designed [30,31]. Two decentralized IMC controllers using the RGA pairing and two decentralized IMC controllers using the RNGA pairing are designed. The notion of a unit step set-point change is applied sequentially to both control loops and adopted the Integral Absolute Error (IAE) [32] performance criterion to examine the efficacy of controllers.

A set of two IMC tuned PID controller parameters that obeys the 1-1/2-2 pairing (RNGA) are listed in Table 2. A similar procedure is adapted to calculate two IMC tuned PID controller parameters that obey the 1-3/2-4 pairing (RGA). The controller structure is

$$G_c(s) = k_c(1 + \frac{1}{\tau_i s} + \tau_d s).$$

**Table 2.** Decentralized PID controllers for both control configurations of radiator example

| Loop | 1-1/2-2 pairing (RNGA) | | | Loop | 1-3/2-4 pairing (RGA) | | |
|---|---|---|---|---|---|---|---|
| | $k_c$ | $\tau_i$ | $\tau_d$ | | $k_c$ | $\tau_i$ | $\tau_d$ |
| 1-1 | -2.434 | 49.305 | 5.912 | 1-3 | 2.697 | 82.576 | 8.279 |
| 2-2 | 1.928 | 38.544 | 6.501 | 2-4 | 2.372 | 68.396 | 7.914 |

Now, using the 1-1/2-2 pairing recommended by RNGA of (15), embedding the IMC tuned PID controller and changing the unit step set-point of the first output $Y_1$, the IAE values of two closed-loops are calculated. Note that the second output is set-point change-free. The IAE of the RNGA-based first output $Y_1$ is less than the IAE of the RGA-based first output. The IAE of the RNGA-based second output $Y_2$ is also less than the IAE of the RGA-based second output, see Figures 5(a) and 5(b). Figures 5(a) and 5(b) show the IAE values, i.e. $(26.67 \quad 51.91)^T$ and $(9.87 \quad 16.15)^T$, respectively.

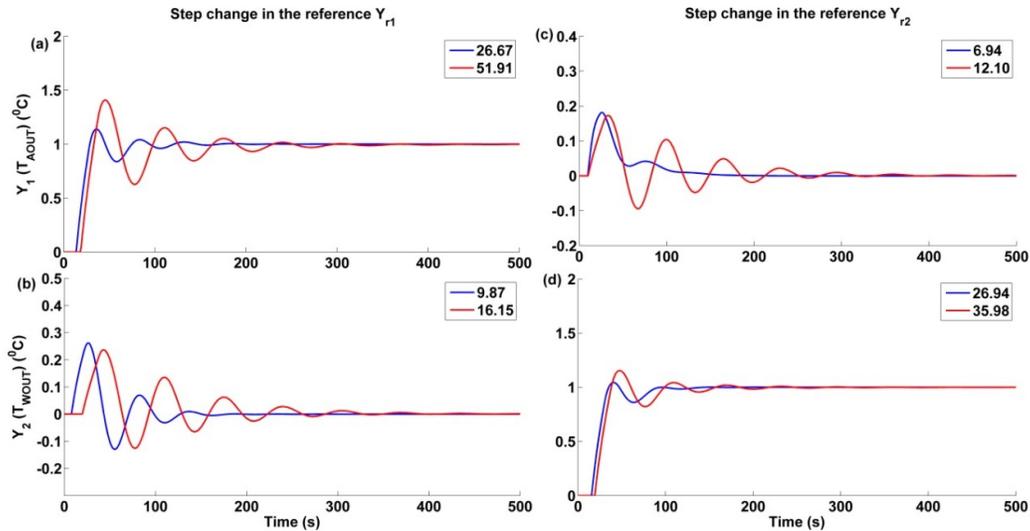

**Figure 5.** Simulation results for radiator example (Blue lines: RNGA recommended pairing, Red lines: RGA recommended pairing)

Figure 5(a) shows the closed-loop response of air outlet temperature $Y_1(T_{AOUT})$ for a step change in the reference $Y_{r1}$. Whereas, Figure 5(b) shows the behavioral pattern of water outlet temperature $Y_2(T_{WOUT})$ under interactions from the change in the reference $Y_{r1}$. Now, consider the case were a unit step set-point change to the second output is given, i.e. change in reference $Y_{r2}$, keeping the set-point of first output unchanged. Figures 5(c) and 5(d) display the same revelation about the IAE values. Figures 5(c) and 5(d) show the IAEs values, i.e. $(6.94 \quad 12.10)^T$ and $(26.94 \quad 35.98)^T$, respectively. Figures 5(a) and 5(d) demonstrates that the closed-loop response resulting from the RNGA-based pairing is better in comparison to the RGA based pairing. It is observed from Figures 5(b) and 5(c) that interactions occurring in the closed-loop performance based on RNGA pairing are less and settles down fast. On the other hand, the interactions in the RGA-based closed-loop are more oscillating.

To evaluate the control efforts put by the controllers to achieve the desired closed-loop performance, Integral Square Control Input (ISCI) is adopted, i.e. $ISCI = \int_o^t u^2(\tau)d\tau$. Less value of ISCI indicates less efforts required by the controller to achieve the desired output [33]. To demonstrate the reduction in control efforts a step set-point change is introduced in the reference $Y_{r1}$. The ISCI index for both RNGA suggested pairing as well as RGA suggested pairing is displayed in Table 3.

**Table 3.** ISCI index for both control configurations of radiator example

|  |  | Control effort in inputs (ISCI) | |
|---|---|---|---|
|  |  | 1-1/2-2 pairing (RNGA) Proposed method | 1-3/2-4 pairing (RGA) |
| Step change in | $Y_1 - Y_{r1}$ | 768.1 | 1092 |
| $Y_{r1}$ | $Y_2 - Y_{r1}$ | 30.3 | 49.97 |

The values depicted in Table 3 shows better performance and less efforts required from the control configuration selected through RNGA suggested pairing in lieu of the RGA suggested control configuration. Hence, the IAE values, Figure 5 and ISCI values are indicative of the superiority of the IMC tuned PID controllers obeying the RNGA-based control-loop pairing in contrast to the RGA based control-loop pairing. Thus, the proposed non-square RNGA of the paper gives a better suggestion of the control-loop pairing for minimum interactions amongst the loop with less efforts required by the controller.

**Conclusion**

In this paper, a formal theory of the non-square RNGA for the 'less output-more input multivariable systems' with the systematic derivation of proofs of the non-square RNGA $\lambda^{RN}$ properties for $r < s$ case is proposed. The theory of RNGA developed in this paper is successfully applied to a non-square radiator laboratory setup with four inputs and two outputs. The results of the experiment carried out on the radiator setup are utilized to obtain the non-square transfer matrix. Then by applying the theory of RNGA to this non-square transfer matrix the control configuration for the decentralized control is achieved. This proves the usefulness of the proposed method to real field practical problems for decentralized and decoupling control. The closed-loop performance of the radiator setup resulting from the RNGA-based controller configurations is compared to that of the RGA based controller configuration. The numerical simulation results of this paper reveal that the controllers of the non-square radiator multivariable system resulting from the RNGA-based pairing are superior to the RGA-based. This indicates the usefulness of RNGA based pairing over RGA based pairing for minimum interactions and better control of non-square multivariable systems. The method proposed in this paper is suggestive for field engineers dealing with the control problem of non-square multivariable systems.

**Acknowledgment**

The Authors express their gratefulness to one of their Colleagues, Dr. Manish Rathod of the Mechanical Engineering Department of the Institute, for useful discussions and providing support for the radiator experimentation in the laboratory environment.